\documentclass[11pt]{amsart}
\usepackage{hyperref}
\usepackage[active]{srcltx}

\usepackage{amsmath, amsfonts,amsthm,times,graphics,color}
\newcommand{\vertiii}[1]{{\left\vert\kern-0.25ex\left\vert\kern-0.25ex\left\vert #1
    \right\vert\kern-0.25ex\right\vert\kern-0.25ex\right\vert}}
 \makeatletter
\renewcommand*\subjclass[2][2000]{%
  \def\@subjclass{#2}%
  \@ifundefined{subjclassname@#1}{%
    \ClassWarning{\@classname}{Unknown edition (#1) of Mathematics
      Subject Classification; using '1991'.}%
  }{%
    \@xp\let\@xp\subjclassname\csname subjclassname@#1\endcsname
  }%
}
 \makeatother
\usepackage{enumerate,amssymb,  mathrsfs,yhmath}

\newtheorem{theorem}{Theorem}[section]
\newtheorem{lemma}[theorem]{Lemma}
\newtheorem*{lemma*}{Lemma}
\newtheorem{proposition}[theorem]{Proposition}
\newtheorem{corollary}[theorem]{Corollary}

\def\1ton{1,2,\ldots,n}
\def\det{{\rm det}}

\usepackage{amssymb}
\usepackage{amsthm}
\usepackage{mathrsfs, amsfonts, amsmath}
\usepackage{graphicx}

\newcommand{\csch}{\operatorname{csch}}

\newtheorem{conjecture}[theorem]{Conjecture}

\newtheorem{remark}[theorem]{Remark}

\numberwithin{equation}{section}

\newcommand{\abs}[1]{\lvert#1\rvert}

\def\XXint#1#2#3{{\setbox0=\hbox{$#1{#2#3}{\int}$}
\vcenter{\hbox{$#2#3$}}\kern-.5\wd0}}

\def\ge{\geqslant}
\begin{document}

\title[Isoperimetric inequality for nearly spherical domains in the Bergman ball]{Isoperimetric inequality for nearly spherical domains in the Bergman ball in $\mathbb{C}^n$}


\author{David Kalaj}
\address{University of Montenegro, Faculty of Natural Sciences and
Mathematics, Cetinjski put b.b. 81000 Podgorica, Montenegro}
\email{davidk@ucg.ac.me}


\begin{abstract}
 We prove a quantitative isoperimetric inequality for the nearly spherical subset of the Bergman ball in $\mathbb{C}^n$. We prove the Fuglede theorem for such sets. This result is a counterpart of a similar result obtained for the hyperbolic unit ball and it makes the first result on the isoperimetric phenomenon in the Bergman ball.
\end{abstract}
\maketitle
\tableofcontents
\sloppy

\maketitle
\section{Introduction and the formula for the perimeter}

The \textit{isoperimetric inequality} states that among all subsets of \( \mathbb{R}^m \) with a given volume, the sphere minimizes the perimeter. More precisely,
\[
P(E) \geq C_m |E|^{\frac{m-1}{m}},
\]
where \( P(E) \) is the perimeter, \( |E| \) is the volume, and
\[
C_m = m \, \omega_m^{1/m},
\]
with $\omega_m$ denoting the volume of the unit ball in $\mathbb{R}^m$. Equality holds if and only if \( E \) is a sphere.

\subsection*{Fuglede's Stability Theorem}\cite{fug}
While the classical inequality identifies spheres as optimal shapes, \textit{Fuglede's theorem} addresses the question of stability:
\textit{How close is a set \( E \) to a sphere if its perimeter nearly minimizes the isoperimetric inequality?}

Fuglede's result provides a quantitative estimate linking the deficit
\[
\delta(E) = P(E) - C_m |E|^{\frac{m-1}{m}},
\]
to the geometric deviation of \( E\) from a sphere. Specifically,
\begin{itemize}
    \item if \( \delta(E) \) is small, then \(E \) is close (in an appropriate metric) to a spherical domain with the same volume.
\end{itemize}

The paper \cite{nicola} by Fusco, Maggi, and Pratelli is a substantial contribution to the study of isoperimetric inequalities. It advances the understanding of how the isoperimetric deficit, defined as
\[
D(E) = \frac{P(E)}{C_m |E|^{(m-1)/m}} - 1,
\]
and the Fraenkel asymmetry,
\[
\lambda(E) = \min \left\{ \frac{|E \Delta (x + rB)|}{r^m} : x \in \mathbb{R}^m, \ |rB| = |E| \right\},
\]
are related. The results are pivotal in the field and significantly refine the classical isoperimetric inequality.

\subsection*{Extensions to hyperbolic space}
The extension of the isoperimetric inequality to the hyperbolic space has been given by Schmidt in \cite{14}.
In the \textit{hyperbolic space \( \mathbb{H}^m \)}, the geometry is non-Euclidean, and the isoperimetric inequality takes a different form, reflecting the curvature of the space. Here, geodesic spheres minimize the isoperimetric ratio.

V. B\"ogelein, F. Duzaar, C. Scheven \cite{duzzar} extended stability results to hyperbolic settings. They showed that:
\begin{enumerate}
    \item a quantitative stability result analogous to Fuglede's holds in \( \mathbb{H}^m \);
    \item the deficit in the hyperbolic isoperimetric inequality controls the deviation of a domain from a geodesic sphere.
\end{enumerate}

The Bergman ball in \(\mathbb{C}^n\) serves as a counterpart to hyperbolic space. In the case \(n = 1\) (i.e., the complex unit disk), these two spaces coincide, but for \(n > 1\), they exhibit significant differences.
A natural and important approach to the geometric analysis of the Bergman ball is therefore to establish an isoperimetric inequality for subsets of this space, namely for the unit ball endowed with the Bergman metric.
Determining such an inequality would, in particular, lead to the computation of the isoperimetric profile of the Bergman ball (for a discussion of isoperimetric profiles, see \cite{ritore}).

At present, however, a complete isoperimetric inequality in the Bergman ball remains a challenging open problem. The main difficulty lies in the intrinsic structure of the Bergman metric, whose anisotropic and non-conformally flat nature makes many classical geometric techniques unavailable.
A particularly interesting question is whether the extremal sets for such an inequality are Euclidean balls centered at the origin and their holomorphic images.
Evidence supporting this conjecture is provided in this paper, but we do not prove the full isoperimetric inequality.
Instead, our contribution is a \emph{Fuglede-type quantitative stability result} in this setting, showing that if a set nearly saturates the expected inequality, then it must be close (in a suitable sense) to the conjectured extremal configuration.

A recent attempt at solving the full isoperimetric problem appears in \cite{suli}. Moreover, proving such an isoperimetric inequality would have implications for the Lieb--Wehrl entropy conjecture in the context of \(SU(N,1)\) (see \cite{prame} and also \cite{NicolaRiccardiTilliWehrl}). In the one-dimensional case, the result was established by Kulikov in \cite{kulik} using the isoperimetric inequality for the hyperbolic disk. For some extension of Kulikov result see \cite{kalaj-jfa} and \cite{kalaj-nyjm}.

The stability framework unifies understanding across Euclidean and non-Euclidean geometries, emphasizing the robustness of optimal shapes under perturbations. The Faber--Krahn inequality is essentially a spectral isoperimetric inequality. It emerges from the same principle as the isoperimetric inequality: balls are optimal shapes for fundamental optimization problems in geometry and analysis. The stability of the Faber--Krahn inequality for the short-time {Fourier} transform published in \cite{inventnew} is also based on the isoperimetric inequality for the Euclidean space. For an extension and improvement se the paper  \cite{inventnew}, see \cite{jkmr}.

\subsection{Bergman unit ball}
Let $\mathbb{B}_n=\{z\in\mathbb{C}^n: \vert z\vert<1\}$ be the unit ball in $\mathbb{C}^n$.
The boundary of $\mathbb{B}_n$ will be denoted by $\mathbb{S}$. Let $ \upsilon$ denote the Euclidean volume measure on $\mathbb{B}_n$.
The $2n-1$ dimensional Hausdorff measure will be denoted by $\mathcal{H}$ and by $d\mathcal{H}$ its area element.  By $\mathbb{B}_r$ we denote the ball $|z|<\tanh(r/2)$.

The Bergman metric in the unit complex ball \( \mathbb{B}_n = \{ z \in \mathbb{C}^n : | z | < 1 \} \) is defined as follows (\cite{rudin}):
 First, define the Bergman kernel on \( \mathbb{B}_n \) by
   \[
   K(z, w) =  \frac{1}{(1 - \langle z, w \rangle)^{n+1}}
   \]
   where \( \langle z, w \rangle = z_1 \overline{w_1} + z_2 \overline{w_2} + \cdots + z_n \overline{w_n} \) is the standard Hermitian inner product on \( \mathbb{C}^n \). The real inner product is then defined  by \( \langle z, w \rangle_R =\Re  \langle z, w \rangle\).
The Bergman metric is defined by the matrix \( B(z) = \{ g_{i\bar{j}}(z) \}_{i,j=1}^n \), where
\[
g_{i \bar{j}}(z) = \frac{\partial^2}{\partial z_i \, \partial \overline{z_j}} \log K(z, z),
\]
which provides the components of the metric tensor.

In the unit ball \( \mathbb{B}_n \), this can be computed explicitly as:
   \[
   g_{i \bar{j}}(z) = \frac{(1 - | z |^2) \delta_{ij} + \overline{z_i} z_j}{(1 - | z |^2)^2}.
   \]
Its inverse is given by the following entries,
   \[
   g^{i \bar{j}}(z) = (1 - | z |^2) (\delta_{ij} - \overline{z_i} z_j).
   \]
\subsection{The Bergman perimeter of a set in the Bergman ball} The aim of this section is to find the formula for the perimeter of a level set in the Bergman ball.
The measure (fundamental form) on the Bergman ball is defined by $$d\mu(z) = (1-|z|^2)^{-n-1} d\upsilon(z),$$ where $d\upsilon(z)$ is the $2n$ dimensional Lebesgue measure. So the Bergman measure of an set $E\subset \mathbb{B}$ is given by $$\mu(E) = \int_{E} d\mu (z).$$
   Assume that $U$ is a smooth real function defined in the unit ball and assume that $$M=\{z\in\mathbb{B}: U(z)=c\}$$ is compactly supported in $\mathbb{B}$.

  For a point $z\in M$ let $N=\nabla_b U$ be the induced normal by the Bergman metric at $M$. Assume $\sigma_g$ is the volume form in $\mathbb{B}$ and let $\sigma_{\tilde g}$ be the induced area form in $M$. We use the formula
\begin{equation}\label{tilde}
\tilde \sigma=\sigma_{\tilde g}=\iota^\ast_M(N\lrcorner \sigma_g)\end{equation} see \cite[Proposition~15.32]{lee}.
  Now we find the unit normal vector field. We first have  $$ \nabla_b U= \sum_{i=1}^n g^{i\bar{j}}(z) \frac{\partial U}{\partial z_i}\frac{\partial}{\partial z_j}.$$
Then
$$\nabla_b U = (1-|z|^2)((1-|z_1|^2)U_{\overline{z}_1}-\sum_{k=2}^n z_1\overline{z_k}U_{\overline{z}_k},\dots,
(1-|z_n|^2)U_{\overline{z}_n}-\sum_{k=1}^{n-1} z_n\overline{z_k}U_{\overline{z}_k}).$$
   We compute
\[\begin{split}\left<\nabla_b U, \nabla U\right> &= (1-|z|^2)\left((1-|z|^2) |\nabla U|^2 + \sum_{1\le j<k\le n}|z_k U_{\bar{z}_j}-z_j U_{\bar z_k}|^2\right)
\\&=(1-|z|^2)\left( (1-|z|^2) |\nabla U|^2 + |z|^2 |\nabla U|^2 - |\left< z,\nabla U\right>|^2\right)
\\&=  (1-|z|^2)(|\nabla U|^2-|\left< z,\nabla U\right>|^2). \end{split}\]
We also have $$|\nabla_b U|_b^2 = \left<\nabla_b U, \nabla U\right>.$$

This implies the formula \begin{equation}\label{formu}
\left<\frac{\nabla_b U}{|\nabla_b U|_b}, \frac{\nabla U}{|\nabla U|}\right>= \sqrt{1-|z|^2}\sqrt{1-\left|\left< \frac{\nabla U}{|\nabla U|}, z\right>\right|^2}.
\end{equation}
The inner product between two vectors $V_1, V_2$ in the Bergman ball is given by $$\left<V_1, V_2\right>_b= \Re[\left< B(z)V_1, V_2\right>].$$ Now for $V_1=\nabla_b U$, we have $B(z) V_1 = c \nabla U,$ with $c>0$. From this we see that, $V\in T_zM$ w.r.t. Bergman structure if and only if  $V\in T_zM$ w.r.t. Euclidean structure of $\mathbb{B}_n\subset \mathbb{C}^n$.

Then from \eqref{tilde} and \eqref{formu}, for vectors $v_1,\dots, v_{2n-1}$ from $T_zM$ we find that
\[\begin{split} \tilde \sigma(v_1, \dots, v_{2n-1})&=\mu\left(\frac{\nabla_b U}{|\nabla_b U|_b}, v_1, \dots, v_{2n-1}\right)
\\& =(1-|z|^2)^{-n-1}\upsilon\left(\frac{\nabla_b U}{|\nabla_b U|_b}, v_1, \dots, v_{n-1}\right)
 \\&= (1-|z|^2)^{-n-1}\left<\frac{\nabla_b U}{|\nabla_b U|_b}, \frac{\nabla U}{|\nabla U|}\right> \upsilon\left(\frac{\nabla U}{|\nabla U|},v_1, \dots, v_{2n-1}\right)
\\&=(1-|z|^2)^{-n-1} \left<\frac{\nabla_b U}{|\nabla_b U|_b}, \frac{\nabla U}{|\nabla U|}\right> \sigma(v_1, \dots, v_{2n-1})
\\&= (1-|z|^2)^{-n-\frac{1}{2}} \sqrt{1-\left|\left< \frac{\nabla U}{|\nabla U|}, z\right>\right|^2} \sigma(v_1, \dots, v_{2n-1}).
\end{split}\]

Here $\upsilon$ is Euclidean volume, and $\mu$ is the volume associated to the Bergman metric and  $\sigma(v_1, \dots, v_{2n-1})$ is the $2n-1$ Euclidean volume of the parallelepiped generated by $\{v_1,\dots, v_{2n-1}\}$.

Then we get the following formula which is crucial for this paper \begin{equation}\label{forperimeter}P(E) = \int_{U(z) =c}(1-|z|^2)^{-n-\frac{1}{2}} \sqrt{1-\left|\left< \frac{\nabla U}{|\nabla U|}, z\right>\right|^2}d\mathcal{H}(z).\end{equation} Remark, that a special case of this formula for the balls $\mathbb{B}_r=\{z: |z|<t=\tanh(r/2)\}$, concentric to the Bergman ball has been calculated by Kehe Zhu in \cite[Sect.~1.6]{kezu}. His formula reads as follows
\begin{equation}\label{br}
P(\mathbb{B}_r) = \int_{\mathbb{S}_r}t^{2n-1}(1-t^2)^{-n} d\mathcal{H}(\omega),
\end{equation}
where $\mathbb{S}_r=\{z: |z|=r\}$, which coincides with \eqref{forperimeter} for $U(z) = |z|=\tanh(r/2)$, because in this particular case $$\left< \frac{\nabla U}{|\nabla U|}, z\right>=\left<\frac{z}{|z|},z\right>=|z|.$$

For the case \( n=1 \), the above formula coincides with the perimeter of a curve \( \gamma \) in the hyperbolic unit disk, a well-known result (see \cite{14} and its application in \cite{kulik}). Specifically, in this case,
\[
\left\langle \frac{\nabla U}{|\nabla U|}, z \right\rangle = \frac{U_{\bar{z}}}{|U_{\bar{z}}|} \overline{z},
\]
and thus,
\[
\left| \left\langle \frac{\nabla U}{|\nabla U|}, z \right\rangle \right|^2 = |z|^2.
\] Then $$P(\gamma)=\int_{\gamma}\frac{|dz|}{(1-|z|^2)}.$$

\subsection{Automorphisms of the unit ball $B\subset \mathbb{C}^n$}\label{sub2}

Let $P_a$ be the
orthogonal projection of $\mathbb{C}^n$ onto the subspace $[a]$ generated by $a$, and let $$Q=Q_a =
I - P_a$$ be the projection onto the orthogonal complement of $[a]$. Explicitly, $P_0 = 0$ and  $P=P_a(z) =\frac{\left<z,a\right> a}{\left<a, a\right>}$. Set $s_a = (1 - |a|^2)^{1/2}$ and consider the map
\begin{equation}
\label{Bergman_transf}
p_a(z) =\frac{a-P_a z-s_a Q_a z}{1-\left<z,a\right>}.
\end{equation}
Compositions of mappings of the form (\ref{Bergman_transf}) and unitary linear mappings of the $\mathbb{C}^n$ constitute the group of holomorphic automorphisms of the unit ball $\mathbb{B}_n\subset \mathbb{C}^n$. It is easy to verify that $p_a^{-1}=p_a$. Moreover, for any automorphism $q$ of the Bergman ball onto itself there exists a unitary transformation $U$ such that \begin{equation}\label{automob}p_{q(c)}\circ q=U\circ p_c.\end{equation}

By using the representation formula \cite[Proposition~1.21]{kezu}, we can introduce the Bergman metric as \begin{equation}\label{bergmet}d_b(z,w)=\frac{1}{2}\log\frac{1+|p_w(z)|}{1-|p_w(z)|}.\end{equation}
If $\Omega= \{z\in \mathbb{C}^n:\left<z,a\right>\neq 1\}$,  then the map $p_a$  is holomorphic in $\Omega$. It is clear that $\overline{\mathbb{B}_n}\subset \Omega$ for $|a| < 1.$

It is well-known that every automorphism $q$ of the unit ball is an isometry with respect  to the Bergman metric, that is: $d_b(z,w)=d_b(q(z),q(w))$. In particular, every isometry preserves the volume and perimeter of any measurable subset $E\subset \mathbb{B}$ of finite measure.
 \subsection{Holomorphic barycenter of a set}

We know from \cite[p.~54]{goldman} that the Bergman metric has negative (nonconstant) sectional curvature, which implies that the square of the corresponding distance function is strictly geodesically convex (see \cite[eq.~4.8.7]{jost}). Consequently, the mapping
\begin{equation}\label{baryb}
    K(a) = \int_{E} \log \cosh^2 (d_b(z,a)) (1 - |z|^2)^{-n-1} \, d\upsilon(z)
\end{equation}
is strictly geodesically convex. Therefore, it attains a unique minimum, which we call  the \emph{holomorphic barycenter} of $E$.

Let $c$ be the barycenter of $E$. Then, the gradient of $K(a)$ vanishes at $a = c$, leading to the equation
\begin{equation}\label{barizero}
    \int_{E} p_c(z) (1 - |z|^2)^{-n-1} \, d\upsilon(z) = 0.
\end{equation}
 The barycenter is \emph{holomorphically invariant}; in other words, if $\phi$ is any biholomorphic map of the unit ball onto itself and $a$ is the barycenter of a set $E \subset \mathbb{B}^n$, then $\phi(a)$ is the barycenter of the image $\phi(E)$. This reflects the natural compatibility of the barycenter with the holomorphic automorphism group of the ball. For further details, we refer the reader to the recent paper \cite{jk}.

Alternatively, one could construct a barycenter using the strict convexity of the function
\[
    H(c) = \int_{E} d^2_b(z,c) (1 - |z|^2)^{-n-1} \, d\upsilon(z),
\]
as in \cite{duzzar}. However, we choose to maintain our definition since the corresponding stationary equation takes a more convenient form.

Including this introduction, the paper consists of four additional sections. In Section 2, we formulate the main results. Section 3 presents auxiliary results necessary for the proofs, which hold for an arbitrary complex dimension \( n \). In Section~\ref{ksec}, we establish a key inequality for the operator \( \nabla_{\imath z} \), which appears in the formula for the perimeter. Finally, Section 5 contains the proof of Theorem~\ref{mainth}.

\section{The main result}
To formulate our main result, let us introduce some notation. Let $\mathbb{S}$ be the unit sphere in $\mathbb{C}^n$.
For a given point $z \in \mathbb{S}$, the tangent space $T_z \mathbb{S}$ is the set of all vectors in $\mathbb{C}^n$ that are orthogonal to $z$, i.e.,
\[
T_z \mathbb{S} = \{ v \in \mathbb{C}^n : \operatorname{Re}(\langle z, v \rangle) = 0 \},
\]
where $\langle z, v \rangle = \sum_{j=1}^n z_j \overline{v_j}$ is the Hermitian inner product.

Let $\tau = \{\tau_1, \dots, \tau_{2n-1}\}$ be an orthonormal frame for $T_z \mathbb{S}$, consisting of vectors that satisfy:
\[
\langle \tau_i, \tau_j \rangle = \delta_{ij}, \quad \text{and } \operatorname{Re}(\langle z, \tau_i \rangle) = 0 \text{ for all } i.
\]

Let $\nabla$ denote the Levi-Civita connection on the sphere $\mathbb{S}$. For any real valued differentiable function $U$ on $\mathbb{S}$, the covariant derivative in the direction of $\tau_i$ is defined as $\nabla_{\tau_i} U$. The notation $\nabla_\tau U$ refers to the collection of these covariant derivatives:
\[
\nabla_\tau U = \{ \nabla_{\tau_1} U, \nabla_{\tau_2} U, \dots, \nabla_{\tau_{2n-1}} U \}.
\]

The norm of $\nabla_\tau$ is given by the sum of the squared norms of these directional derivatives:
\[
|\nabla_\tau U|^2 = \sum_{i=1}^{2n-1} (\nabla_{\tau_i} U)^2.
\]
Then we define $$\|U\|^2_{W^{1,2}}=\int_{\mathbb{S}}|\nabla_\tau U|^2 d\mathcal{H} + \int_{\mathbb{S}}|U|^2 d\mathcal{H}.$$ We also define $$\|U\|_{W^{1,\infty}}=\max\{|\nabla_\tau U|_\infty, |U|_\infty\}.$$

\begin{theorem}[Fuglede's theorem in the Bergman ball]\label{mainth}
For every $r_0>0$ there exists $\varepsilon_0\in\bigl(0,\tfrac12\bigr]$, depending only on $r_0$,
such that the following holds.

Let $E\subset \mathbb{B}_n$ be a set whose barycenter is at the origin and such that
\begin{equation}\label{eq:vol_constraint_thm}
\mu(E)=\mu(\mathbb{B}_r)
\end{equation}
for some $r\in(0,r_0]$. Assume that $\partial E$ is a radial graph over the unit sphere
$\mathbb{S}$, namely there exists a Lipschitz function $u:\mathbb{S}\to\mathbb{R}$ with
\begin{equation}\label{eq:graph_param_thm}
Z(\omega)=\omega\,\tanh\!\left(\frac r2\bigl(1+u(\omega)\bigr)\right),\qquad \omega\in\mathbb{S},
\end{equation}
such that $Z(\mathbb{S})=\partial E$, and
\begin{equation}\label{eq:u_small_thm}
\|u\|_{W^{1,\infty}(\mathbb{S})}\le \varepsilon_0.
\end{equation}
Then there exists a constant $c_1=c_1(r_0,n)>0$ such that
\begin{equation}\label{eq:fuglede_estimate}
\frac{P(E)-P(\mathbb{B}_r)}{P(\mathbb{B}_r)}
\;\ge\;
c_1\,\|u\|_{W^{1,2}(\mathbb{S})}^2.
\end{equation} In particular,
\[
P(E)\ge P(\mathbb B_r),
\]
with equality if and only if $E=\mathbb B_r$.
\end{theorem}

\begin{corollary}[Stability for nearly spherical ellipsoids in the Bergman ball]\label{thm:ellipsoid_stability}
Fix $r_0>0$. Let $\varepsilon_0\in\bigl(0,\frac12\bigr]$ and $c_1(r_0)>0$ be the constants
from Theorem~\ref{mainth}. Then there exists $\delta_0=\delta_0(n,r_0)>0$ such that the
following holds.

Let $E\subset\mathbb B_n$ be a Euclidean ellipsoid centered at the origin of the form
\[
E=\left\{z\in\mathbb C^n:\ \sum_{j=1}^n \frac{|z_j|^2}{\lambda_j^2}<R^2\right\},
\qquad \lambda_j>0,
\]
with
\[
\max_{1\le j\le n}|\lambda_j-1|\le \delta_0.
\]
Assume that $R$ is chosen so that
\[
\mu(E)=\mu(\mathbb B_r)
\]
for some $r\in(0,r_0]$, and assume moreover that $E$ has barycenter at the origin.

Then $\partial E$ can be written as a nearly spherical graph over $\mathbb S$ in the sense
of Theorem~\ref{mainth}, and the quantitative isoperimetric estimate
\[
\frac{P(E)-P(\mathbb B_r)}{P(\mathbb B_r)}\ge c_1(r_0)\,\|u\|_{W^{1,2}(\mathbb S)}^2
\]
holds. 
\end{corollary}

\begin{proof}
Since $E$ is centered at the origin and strictly star-shaped, its boundary can be written
as a radial graph over the unit sphere:
\[
\partial E=\{\rho(\omega)\,\omega:\ \omega\in\mathbb S\}.
\]
Writing $\omega=(\omega_1,\dots,\omega_n)\in\mathbb C^n$ with $|\omega|=1$, the equation
$\rho(\omega)\omega\in\partial E$ gives
\[
\sum_{j=1}^n \frac{|\rho(\omega)\omega_j|^2}{\lambda_j^2}=R^2,
\]
hence
\[
\rho(\omega)=\frac{R}{\sqrt{\sum_{j=1}^n \frac{|\omega_j|^2}{\lambda_j^2}}}.
\]
Assume $\max_j|\lambda_j-1|\le \delta$ with $\delta\in(0,1)$ small. Then
\[
\lambda_j^{-2}=1+O(\delta)
\quad\text{uniformly in }j,
\]
and using $\sum_{j=1}^n|\omega_j|^2=1$ we infer
\[
\sum_{j=1}^n\frac{|\omega_j|^2}{\lambda_j^2}=1+O(\delta)
\quad\text{uniformly on }\mathbb S,
\]
so that
\[
\rho(\omega)=R\,(1+O(\delta))
\quad\text{uniformly on }\mathbb S.
\]
Differentiating the formula for $\rho$ on $\mathbb S$ yields in addition
\[
\|\nabla_{\mathbb S}\rho\|_{L^\infty(\mathbb S)}\le C(n)\,R\,\delta.
\]
Consequently,
\[
\|\rho-R\|_{W^{1,\infty}(\mathbb S)}\le C(n)\,R\,\delta.
\]
Choosing $\delta_0>0$ so that $C(n)\delta_0\le \varepsilon_0$ implies that $\partial E$
admits a parametrization of the form
\[
Z(\omega)=\omega\,\tanh\!\left(\frac r2(1+u(\omega))\right),
\]
for some $u\in W^{1,\infty}(\mathbb S)$ satisfying
\[
\|u\|_{W^{1,\infty}(\mathbb S)}\le \varepsilon_0.
\]
Thus $E$ is nearly spherical in the sense of Theorem~\ref{mainth}.

\medskip

Since $\mu(E)=\mu(\mathbb B_r)$ and the barycenter of $E$ is at the origin, all
assumptions of Theorem~\ref{mainth} are satisfied. Therefore,
\[
\frac{P(E)-P(\mathbb B_r)}{P(\mathbb B_r)}\ge c_1(r_0)\,\|u\|_{W^{1,2}(\mathbb S)}^2\ge 0,
\]
which yields $P(E)\ge P(\mathbb B_r)$.
If equality holds, then $\|u\|_{W^{1,2}(\mathbb S)}=0$, hence $u$ is constant, and the
volume constraint forces $u\equiv 0$, which implies $E=\mathbb B_r$.
\end{proof}

\begin{remark}
We can assume in the previous theorem that $r_0\ge \cosh^{-1}\left[\frac{-2+3 n+\sqrt{1-8 n+8 n^2}}{-1+n}\right]$. For that $r_0$, we can choose $$c_1(r_0)=\frac{r_0^2}{2\omega_{2n-1} \sinh^2r_0 }.$$

After completing this paper, we realized that a different lower bound estimate for 
\(P(E) - P(\mathbb{B}_r)\) was previously obtained by Silini in \cite{Silini2025} 
for the class of rank-one symmetric spaces of non-compact type. While his result is 
more general, our result is \emph{more quantitative} and applies directly to the 
\emph{isoperimetric deficit}. Both results, whose proofs are different, however, lead to the same qualitative 
conclusion: when a set is nearly spherical, it must in fact be spherical if it 
\emph{minimizes the perimeter} among all subsets of the unit ball having the same 
\emph{Bergman volume}.
\end{remark}
\section{Auxiliary results}\label{sect3}

\subsubsection{A formula for the volume of $E$.}\label{subsub} A parametrization of $E$ is given by
\[
\tilde{Z} : \mathbb{S} \times [0, r] \to E, \quad \tilde{Z}(\omega, \rho) = \omega \tanh\left(\frac{\rho}{2} (1 + u(\omega))\right).
\]

For a fixed \(\omega \in \mathbb{S}\) and \(\rho \in [0, r]\), we introduce the notation
\[
\mathbf{t} := \tanh\left(\frac{\rho}{2} (1 + u(\omega))\right).
\]
Assume that \(\{e_1, \dots, e_{2n}\}\) is a standard orthonormal basis of \(\mathbb{R}^{2n}\), and let \(\omega = (\omega_1, \dots, \omega_{2n})\).

To compute the Jacobian of \(\tilde{Z}\), we choose an orthonormal real basis \(\tau_1, \dots, \tau_{2n-1}\) of \(T_{\omega} \mathbb{S}\) with respect to the inner product
\[
\langle w, z \rangle_R := \Re \langle w, z \rangle.
\]
In particular, we set
\[
\tau_{2n-1} = \imath\omega = \imath (\omega_1 + \imath\omega_2, \dots, \omega_{2n-1} + \imath\omega_{2n}) = (-\omega_2, \omega_1, \dots, -\omega_{2n}, \omega_{2n-1}) \in \mathbb{S},
\]
and define \(\tau_{2n} = \omega\), ensuring that \(\langle \omega, \imath\omega \rangle_R = 0\).

Next, we construct an orthonormal basis for \(T_{(\omega, \rho)} (\mathbb{S} \times (0, r))\) as
\[
(\tau_1, 0), \dots, (\tau_{2n-1}, 0), e_{2n}.
\]
Similarly, an orthonormal basis for the tangent space at the point $X(\omega, \rho)$: \(T_{X(\omega, \rho)} \mathbb{B}_n \cong \mathbb{R}^{2n}\) is given by
\[
\tau_1, \dots, \tau_{2n-1}, \omega.
\]

Concerning  these choices, the Jacobi matrix of $\tilde{Z}$ reads
\[
D\tilde{Z}(\omega, \rho) =
\begin{pmatrix}
\mathbf{t} & \cdots  & 0 & 0 \\
0 & \cdots  & 0 & 0 \\
\vdots  & \ddots & \vdots & \vdots \\
0 & \cdots  & \mathbf{t} & 0 \\
\frac{\rho(1-\mathbf{t}^2) \partial_1 u}{2}  & \cdots & \frac{\rho(1-\mathbf{t}^2) \partial_{2n-1} u}{2} & \frac{1+u}{2}(1-\mathbf{t}^2)
\end{pmatrix} \in \mathbb{R}^{2n \times 2n}.
\]

Consequently, the Jacobian of $\tilde{Z}$ is given by
\begin{equation}\label{ztil}
J\tilde{Z}(\omega, \rho) = \frac{1 + u}{2} \mathbf{t}^{2n-1} (1 - \mathbf{t}^2).
\end{equation}
The volume of $E$ can thus be written as
\[
\mu(E) = \int_E \left(1 - |z|^2\right)^{-n-1} d\upsilon(z) = \int_{\mathbb{S}} \int_0^{r} \left(\frac{1}{1 - \mathbf{t}^2}\right)^n \frac{1 + u}{2} \mathbf{t}^{2n-1}  \, d\rho \, d\mathcal{H}(\omega).
\]
Recalling the definition of $\mathbf{t}$, we compute
\[
\left(\frac{1}{1 - \mathbf{t}^2}\right)^n \frac{1 + u(\omega)}{2} \mathbf{t}^{2n-1} = \frac{1+u}{2}\cosh\left[\frac{1}{2} \rho (1+u)\right] \sinh\left[\frac{1}{2} \rho (1+u(\omega))\right]^{2n-1}.
\]
So \[\begin{split}\mu(E) &= \int_{E}\frac{1}{(1-|z|^2)^{n+1}}d\upsilon(z)\\&=\int_0^{r}\int_{\mathbb{S}}\frac{1+u}{2}\cosh\left[\frac{1}{2} \rho (1+u(\omega))\right] \sinh\left[\frac{1}{2} \rho (1+u(\omega))\right]^{2n-1}d\mathcal{H}(\omega)d\rho
\\&=\int_{\mathbb{S}}\frac{\sinh^{2n}\left[\frac{1}{2} r (1+u(\omega))\right]}{2n}d\mathcal{H}(\omega).\end{split}\] We recall here that $d\mathcal{H}=dH^{2n-1}$ is the $2n-1$ Hausdorff  measure.
\subsubsection{Consequences of the volume constraint $\mu(\mathbb{B}_r)=\mu  (E)$}

Let $$F(u):=\frac{\sinh^{2n}\left[\frac{1}{2} r (1+u)\right]}{2n}.$$
Then $$\int_{\mathbb{S}} F(u(\omega))-F(0)d\mathcal{H}(\omega)=0.$$
Further
\[\begin{split}
F(u) - F(0) &= \frac{1}{4} r \sinh^{2n-2} (\frac{r}{2})\sinh(r)  u\\&+\frac{1}{8} r^2 (-1+n+n \cosh r) \sinh^{2n-2}(\frac{r}{2}) u^2+\frac{F'''(u_1)}{6} u^3,
\end{split}\]
where \[\begin{split}F'''(u_1)&=\frac{1}{4} r^3 \cosh\left[\frac{1}{2} r (1+u_1)\right] \left(1+(n-3) n+n^2 \cosh[r (1+u_1)]\right)\\&\times \sinh^{2n-3}\left[\frac{1}{2} r (1+u_1)\right].\end{split}\]
Note that $u=u(\omega)$ and $u_1=u_1(\omega)$. Then   $$\int_{\mathbb{S}}     u+\frac{r(-1+n+n \cosh r)}{2 \sinh(r)}  u^2+\frac{4F'''(u_1)}{3r \sinh(r)\sinh^{2n-2} (\frac{r}{2})} u^3d\mathcal{H}(\omega)=0.$$
Since  $|u_1|\le 1/2$, \[\begin{split} \frac{4F'''(u_1)}{6r \sinh(r)\sinh^{2n-2} (\frac{r}{2})}\le  C_1(n, r_0), \end{split}\]
where
\[\begin{split}C_1(n,r_0)&=\cosh\left[\frac{1}{2} r_0 (1+1/2)\right] \left(1+(n-3) n+n^2 \cosh[r_0 (1+1/2)]\right)\\&\times \max_{r\in[0,r_0]}\frac{r^2  \sinh^{2n-3}\left[\frac{1}{2} r (1+1/2)\right]}{{6 \sinh(r)\sinh^{2n-2} (\frac{r}{2})} }.\end{split}\]
Let $$c_0(n,r)=\frac{r}{2 \sinh(r)}  (-1+n+n \cosh r).$$
It follows that \begin{equation}\label{estU}|a_0|=\abs{\int_{\mathbb{S}}     u d\mathcal{H}(\omega)}\le c_0(n,r)\int_{\mathbb{S}} u^2 d\mathcal{H}(\omega)+ \varepsilon_0  C_1(n,r_0)\int_{\mathbb{S}} u^2 d\mathcal{H}(\omega).\end{equation}
Since $$\int_{\mathbb{S}} u^2 d\mathcal{H}(\omega)= \|u\|_{L^2(\mathbb{S})}^2\le \|u\|_{W^{1,2}}^2$$ we obtain that there is a constant $C_2(n,r_0)$ so that
\begin{equation}\label{rur}|a_0|\le \varepsilon_0 C_2(n,r_0)\|u\|_{W^{1,2}}.\end{equation}
\subsubsection{Consequence of barycenter constraint} We begin by the following \begin{remark}
For any point $a \in \mathbb{B}^n$, there exists a biholomorphism $\phi_a$ mapping $a$ to the origin. So if a domain $E$ has barycenter $a$, then its image $E'=\phi(a)(E)$ has the barycenter  equal to $0$ and has the same Bergman volume and perimeter as $E$.

\end{remark}

If $z=0$ is the barycenter, then it follows from \eqref{baryb} that
$$\int_E z (1-|z|^2)^{-1-n} d\upsilon(z)=0.$$
Thus from \eqref{ztil} we have $$ \int_0^{r} \int_{\mathbb{S}}(1-{\mathbf{t}}^2)^{-n-1} \omega \frac{(1+u(\omega))}{2}{\mathbf{t}}^{2n-1}(1-\mathbf{t}^2)  \mathbf{t}d\rho d\mathcal{H}(\omega)=0.$$
Hence
$$ \int_0^{r} \int_{\mathbb{S}}(1-\mathbf{t}^2)^{-n} \omega \frac{(1+u(\omega))}{2}\mathbf{t}^{2n}  d\rho d\mathcal{H}(\omega)=0.$$
Thus $$\int_0^{r} \int_{\mathbb{S}}\omega(1+u(\omega))\sinh^{2n}\left[\frac{1}{2} \rho  (1+u(\omega ))\right]d\rho d\mathcal{H}(\omega)=0.$$
Let $$\phi(u)=(1+u)\sinh^{2n}\left[\frac{1}{2} \rho  (1+u)\right]. $$
Since $$\int_0^{r} \int_{\mathbb{S}}\omega\sinh^{2n}\left[\frac{1}{2} \rho  \right]d\rho d\mathcal{H}(\omega)=0,$$ and  $$ \phi(u) = \phi(0)+\phi'(0)u + \frac{\phi''(u_0)}{2!} u^2, $$ for some $u_0=u_0(\omega)\in[-1/2,1/2]$,
we get  $$\int_0^{r} \int_{\mathbb{S}} (\phi'(0)u + \frac{\phi''(u_0)}{2!} u^2)\omega d\rho d\mathcal{H}(\omega)=0.$$
Further
$$\phi'(0)=\sinh \left[\frac{\rho}{2}\right]^{-1+2 n} \left(n \rho \cosh \left[\frac{\rho}{2}\right]+\sinh \left[\frac{\rho}{2}\right]\right)$$
and
\[\begin{split}\phi''(u_0)&= \frac{1}{2} n \rho(1+u_0) \sinh \left[\frac{\rho(1+u_0)}{2}\right]^{-2+2 n} \\ & \times(\rho(1+u_0) (-1+n+n \cosh [\rho(1+u_0)])+2 \sinh [\rho(1+u_0)]),\end{split}\]
and this implies that
$$\phi''(u)\le \phi''(1)=n \rho \sinh^{-2+2 n}(\rho)  (\rho(-1+n+n \cosh(2 \rho))+\sinh(2\rho)).$$
Since
$$r \sinh(\frac{r}{2})^{2 n} \int_{\mathbb{S}} u \omega d\mathcal{H}(\omega) = -\int_0^{r} \int_{\mathbb{S}}\frac{\phi''(u_0)}{2!} u^2\omega d\rho d\mathcal{H}(\omega),$$
 we obtain
 \[\begin{split}
 |r \sinh(\frac{r}{2})^{2 n} \int_{\mathbb{S}} u \omega d\mathcal{H}(\omega)|&\le \int_0^{r} \int_{\mathbb{S}}|\frac{\phi''(1)}{2!} u^2|\omega| d\rho d\mathcal{H}(\omega)
  \\&=\frac{n}{2} r^2 \cosh(r) \sinh^{2n-1}(r) \int_{\mathbb{S}} u^2|\omega| d\rho d\mathcal{H}(\omega).\end{split}\]
Thus for $\omega=(\omega_1,\dots,\omega_{2n})\in \mathbb{S}$ we have for every $j$

\begin{equation}\label{estuom1}\left|\int_{\mathbb{S}} u \omega_j d\mathcal{H} \right|\le C_3(n,r_0)\varepsilon_0\|u\|_{W^{1,2}(\mathbb{S})}.\end{equation}
\subsubsection{Spherical harmonics and $W^{1,2}$ norm of $u$}
For any integer \( k \geq 0 \), let us denote by \( \Phi_{k,i}, \, i = 1, \ldots, d_k \), the \textit{spherical harmonics of order \( k \)}, i.e., the restriction to \( \mathbb{S} \) of the homogeneous harmonic polynomials of degree \( k \), normalized such that
\[
\|\Phi_{k,i}\|_{L^2(\mathbb{S})} = 1, \quad \text{for all } k \text{ and } i \in \{1, \ldots, d_k\}.
\] Here \begin{equation}\label{dimd}
d_k=\binom{k+2n-1}{k}-\binom{k+2n-3}{ k-2}.
\end{equation}
First of all
\[
\int_{\mathbb{S}} 1 \, d\mathcal{H} = \Omega_n,
\]
where \( \Omega_n \) is the \( (2n-1) \)-dimensional measure of the unit sphere \( \mathbb{S} \). For \( i = 1, \ldots, n \), for $z=(x_1+iy_1,\dots, x_n+i y_n)$ it follows that
\[
\int_{\mathbb{S}} x_i^2 \, d\mathcal{H} =\int_{\mathbb{S}} y_i^2 \, d\mathcal{H}= \frac{1}{2n} \int_{\mathbb{S}} |z|^2 \, d\mathcal{H} = \frac{\Omega_n}{2n}.
\]
Hence, we have
\[
\Phi_{0,0} = \frac{1}{\sqrt{\Omega_n}} \quad \text{and} \quad \Phi_{1,i} = \frac{x_i}{\sqrt{\Omega_n}}, \quad \Phi_{1,i+n} = \frac{y_i}{\sqrt{\Omega_n}}
\]
respectively.

The functions \( \Phi_{k,i} \) are eigenfunctions of the Laplace-Beltrami operator on \( \mathbb{S} \). Specifically, for all \( k \) and \( i \),
\[
-\Delta_{\mathbb{S}} \Phi_{k,i} = k(k + 2n - 2) \Phi_{k,i}.
\]

Therefore if we write
\[
u = \sum_{k=0}^\infty \sum_{i=1}^{d_k} a_{k,i} \Phi_{k,i}, \quad \text{where} \quad a_{k,i} = \int_{\mathbb{S}} u \Phi_{k,i} d\mathcal{H},
\]
we have
\[
\|u\|_{L^2(\mathbb{S})}^2 = \sum_{k=0}^\infty \sum_{i=1}^{d_k} a_{k,i}^2.
\]

\bigskip

Further

\[
\int_{\mathbb{S}} |\nabla u|^2 \, d\mathcal{H}
= -\int_{\mathbb{S}} u \, \Delta_{\mathbb{S}} u \, d\mathcal{H}^{n-1}
= \sum_{k=1}^\infty \sum_{i=1}^{d_k} k(k + 2n - 2)a_{k,i}^2.
\]

We therefore have
\begin{equation}\label{ww0}
\|\nabla u\|^2_{L^2 (\mathbb{S})} = \sum_{k=0}^\infty \sum_{i=1}^{d_k} \big[k(k + 2n - 2)\big]a_{k,i}^2
\end{equation} and
\begin{equation}\label{ww}
\|u\|_{W_{1,2}}^2 = \sum_{k=0}^\infty \sum_{i=1}^{d_k} \big[k(k + 2n - 2) + 1\big]a_{k,i}^2.
\end{equation}



It follows from \eqref{estuom1} that
\begin{equation}\label{estuom}\left|a_{1,i}\right|\le C_3(n,r_0)\varepsilon_0\|u\|_{W^{1,2}(\mathbb{S})}, i=1,\dots, 2n.\end{equation}

\subsubsection{Formula of perimeter} Let now set  $\mathbf{t}=\tanh(\frac{r}{2}(1+u(\omega)))$.

The Jacobi matrix of $Z$ is given by
\[
D{Z}(\omega) =
\begin{pmatrix}
\mathbf{t} & 0 & \cdots & 0  \\
0 & \mathbf{t} & \cdots & 0  \\
\vdots & \vdots & \ddots & \vdots  \\
0 & 0 & \cdots & \mathbf{t}  \\
\frac{\rho(1-\mathbf{t}^2) \partial_1 u}{2} & \frac{\rho(1-\mathbf{t}^2) \partial_2 u}{2} & \cdots & \frac{\rho(1-\mathbf{t}^2) \partial_{2n-1} u}{2}
\end{pmatrix} \in \mathbb{R}^{2n \times 2n-1}.
\]
For the Jacobian of $Z$, similarly as in \cite{fug} and \cite{duzzar} ( see subsection  below)  we  get
\[
\begin{split}[J Z(\omega)]^2 &= \mathbf{t}^{2(2n-1)} + \mathbf{t}^{2(2n-2)}(1 - \mathbf{t}^2)^2 \frac{\rho^2}{4} r^2 |\nabla_\tau u|^2
\\&= \mathbf{t}^{2(2n-1)} \left[ 1 + \left( \frac{1 - \mathbf{t}^2}{2\mathbf{t}} \right)^2 r^2 |\nabla_\tau u|^2 \right]. \end{split}
\]
Let $E_1=\{z: U(z) = c\}$. Then by using \eqref{forperimeter}
 \begin{equation}\label{perimeter}P(E_1)=\int_{E_1} d \tilde \sigma(z)=\int_{\mathbb{S}} Z_\ast (d \tilde \sigma)\end{equation} where $$d\tilde \sigma(z) = (1-|z|^2)^{-n-\frac{1}{2}}{\sqrt{1-\frac{|\left<\nabla U, z\right>|^2}{|\nabla U|^2}}}  d\sigma(z).$$
 Further \begin{equation}\label{form}Z_\ast (d \tilde \sigma) = \frac{\mathbf{t}^{(2n-1)} }{(1-\mathbf{t}^2)^{n+1/2}}{\sqrt{1-\frac{|\left<\nabla U, z\right>|^2}{|\nabla U|^2}}}\left[ 1 + \left( \frac{1 - \mathbf{t}^2}{2\mathbf{t}} \right)^2 r^2 |\nabla_\tau u|^2 \right]^{1/2}d\mathcal{H}(\omega).\end{equation}
Now, we construct the real function \( U \), whose level set is precisely the boundary \( E_1 = \partial E \). The function \( U: \mathbb{C}^n \to \mathbb{R} \) is defined as:
\[
U(z) = \frac{2 \tanh^{-1}(|z|)}{1 + u\left(\frac{z}{|z|}\right)},
\]
where \( |z| \) is the Euclidean norm of the vector \( z \in \mathbb{C}^n \), and \( u\left(\frac{z}{|z|}\right) \) is a function that depends on the direction of \( z \) (i.e., on the unit vector \( \frac{z}{|z|} \)).

Given this definition of \( U \), we can describe the set \( E \) as the level set of \( U \):

\[
E = \{ z \in \mathbb{C}^n : U(z) = r \} = \{ Z(\omega) : \omega \in \mathbb{S} \},
\]
where \( \mathbb{S} \) denotes the unit sphere in \( \mathbb{C}^n \), and \( Z(\omega) \) is a mapping defined by:

\[
Z(\omega) = \omega \tanh\left( \frac{r}{2} (1+u(\omega)) \right).
\]

Next, we compute the derivative of \( U \) with respect to the tangential direction \( \nabla_\tau \). The result is:

\[
\nabla_\tau U(z) = \frac{-2 \nabla_\tau u}{(1 + u(\omega))^2} \tanh^{-1}(|z|),
\]
where \( \nabla_\tau \) represents the gradient in the tangential direction (which is the direction of motion along the boundary of the set), and \( \nabla_\tau u \) is the tangential derivative of the function \( u \). This expression gives the rate of change of \( U(z) \) along the tangential direction at the point \( z \in \mathbb{C}^n \).

For vectors $\tau_k$ as in Subsection~\ref{subsub}, we have
$$\nabla_{\tau_k} U(z)=\frac{-2\nabla_{\tau_k} u}{(1+u(\omega))^2}\tanh^{-1}(|z|)$$
and for $\mathbf{n}=z/|z|$

$$\nabla_{\mathbf{n}} U(z)=\frac{2}{(1-|z|^2)(1+u(\omega))}.$$  In the sequel we use the complex inner product $\left<\cdot,\cdot\right>$ and the real inner product $\left<\cdot,\cdot\right>_R=\mathrm{Re}\left<\cdot,\cdot\right>$ and the identification \[\begin{split}\mathbb{R}^{2n}&=\{(x_1,y_1,\dots, x_n, y_n):x_j, y_j\in \mathbb{R}\}\\&\cong \{(\imath x_1+y_1,\dots, \imath x_n+ y_n):x_j, y_j\in \mathbb{R}\}=\mathbb{C}^n.\end{split}\]
Then
$$\nabla U = \sum_{j=1}^{2n}\left<\nabla U,\tau_j\right>_R \tau_j=(\imath\left<\nabla U,\tau_1\right>_R+ \left<\nabla U,\tau_2\right>_R,\dots, \imath\left<\nabla U,\tau_{2n-1}\right>_R+ \left<\nabla U,\tau_{2n}\right>_R).$$

Since
$$\frac{z}{|z|}=(0,\dots, 1+0\cdot  \imath)\in \mathbb{C}^{n},$$ we have
$$\left<\nabla U(z), \frac{z}{|z|}\right > =\imath\left<\nabla U,\tau_{2n-1}\right>+ \left<\nabla U,\tau_{2n}\right>.$$
 and \begin{equation}\label{prim}|\left<z, \nabla U\right>|= |z|\sqrt{|\partial_{\mathbf{n}} U|^2+|\partial_ {\tau_{2n-1}} U|^2}.\end{equation} Moreover
\begin{equation}\label{secon}|\nabla U|^2 = |\partial_{\mathbf{n}} U|^2 + |\nabla_\tau U|^2= |\partial_{\mathbf{n}} U|^2 + \sum_{k=1}^{2n-1}|\partial_{\tau_k} U|^2 .\end{equation}
Let us introduce the shorthand notation for the sequel $$\mathbf{t}=\tanh[r/2(1+u(\omega))],$$ $$R=\frac{\left(\frac{1+u(\omega)}{1-\mathbf{t}^2}\right)^2}{(\tanh^{-1}(\mathbf{t}))^2}=\frac{4}{r^2 \left(1-\mathbf{t}^2\right)^2 }$$ and $$x=\|\nabla_\tau u\|^2,  s=|\nabla_{\tau_{2n-1}}|^2.$$
Let  $$h(x)={\sqrt{1+\frac{r^2 \left(1-\mathbf{t}^2\right)^2 x}{4\mathbf{t}^2}}} $$ and

 $$k(x,s)=\sqrt{1-\frac{\mathbf{t}^2 (R+s)}{R+ x}}.$$
a) By using the elementary inequality $$\sqrt{1+y}\ge 1+\frac{y}{2}-\frac{y^2}{8},$$ by making the substitution $\mathbf{t}=\tanh(\frac{r}{2}(1+u))$ we obtain
$$h(x)\ge 1+\frac{1}{2} r^2 \csch^2(r (1+u)) x-\frac{1}{8} \left(r^4 \csch^4(r (1+u))\right) x^2.$$
By using continuity and the fact that $|u|\le \varepsilon_0<1/2$ for $r\le r_0$ we obtain, \begin{equation}\label{esth}h(x) \ge 1+ \frac{1}{2} x r^2 \csch^2 r - \varepsilon_0 C_{4}(n,r_0) x.\end{equation}

b) Now consider $$k(x,s)=\frac{\sqrt{1-\frac{(R+s) \mathbf{t}^2}{R+ x}}}{\sqrt{1-\mathbf{t}^2}}.$$
Then $$k(x,s) = h(x,s)+s^2\mathbf{r_1}$$
where $$ h(x,s)=\frac{\sqrt{1-\frac{R \mathbf{t}^2}{R+x}}}{\sqrt{1-\mathbf{t}^2}}-\frac{\mathbf{t}^2 s}{2 \left(\sqrt{1-\mathbf{t}^2} (R+x) \sqrt{1-\frac{R \mathbf{t}^2}{R+x}}\right)}$$ and $$\mathbf{r_1}=-\frac{\mathbf{t}^4}{8 \sqrt{1-\mathbf{t}^2} (R+x)^2 \left(1-\frac{(R+s_1) \mathbf{t}^2}{R+x}\right)^{3/2}},$$ and $s_1\in [0,1/2]$. Given the definition of $R$ and $\mathbf{t}\le \tanh(r_0)$ we obtain that  $$|\mathbf{r_1}|\le C_4(n,r_0).$$ Further $$h(x,s)= 1+\frac{\mathbf{t}^2 (x-s)}{2 R \left(1-\mathbf{t}^2\right)}-\frac{s \mathbf{t}^2 \left(2-\mathbf{t}^2\right) x}{4 R^2 \left(1-\mathbf{t}^2\right)^2}+ \mathbf{r_2} x^2 $$ where \[\begin{split}\mathbf{r_2}&=\frac{R^2 \mathbf{t}^4 \left(-3 s \mathbf{t}^2+2 R \left(-1+\mathbf{t}^2\right)-2 x\right)}{16 \sqrt{1-\mathbf{t}^2} (R+x)^5 \left(1-\frac{R \mathbf{t}^2}{R+x}\right)^{5/2}}
\\&- \frac{\mathbf{t}^2 \left(R^2 \left(1-\mathbf{t}^2\right)+s x+R (s+x)\right)}{2\sqrt{1-\mathbf{t}^2} \left(R \left(1-\mathbf{t}^2\right)+x\right) (R+x)^3 \sqrt{1-\frac{R \mathbf{t}^2}{R+x}}}.\end{split}\] Since $$R\ge \frac{4}{r_0^2(1-\tanh(r_0)^2)}, $$  $$\mathbf{t}\le \tanh(r_0),\ \ \ \ 0\le s\le x\le 1/2$$ and $\lim_{R\to  +\infty } \mathbf{r_2}=0$. It follows that there is a constant $C_5(n,r_0)$ so that   $$|\mathbf{r_2}|\le C_5(n,r_0).$$

Using the condition $\|u\|_{W^{1,\infty}}\le \varepsilon_0$, we obtain that there is a constant $C_6(n,r_0)$ so that
\begin{equation}\label{kxs}
k(x,s)\ge 1+ \frac{r^2}{8}  \tanh^2(\frac{r}{2})\left(1-\tanh^2(\frac{r}{2})\right)(x-s) - C_6(n,r_0)\varepsilon_0 x.\end{equation}

By using the equation
$${\sqrt{1-\frac{|\left<\nabla U, z\right>|^2}{|\nabla U|^2}}}=\sqrt{1-|z|^2 \frac{|\partial_{\mathbf{n}} U|^2 + |\nabla_{\tau_{2n-1}} U|^2}{|\partial_{\mathbf{n}} U|^2 + |\nabla_\tau U|^2}}$$ and \eqref{esth} and \eqref{kxs},
we have
\begin{lemma}\label{lema22} Under previous notation, 
if $0<r<r_0$ and $\|u\|_{W^{1,\infty}}<\varepsilon_0$ then
\[\begin{split}{\sqrt{1-\frac{|\left<\nabla U, z\right>|^2}{|\nabla U|^2}}}&\left[ 1 + \left( \frac{1 - \mathbf{t}^2}{2\mathbf{t}} \right)^2 r^2 |\nabla_\tau u|^2 \right]^{1/2}
\\&\ge 1+\frac{1}{2}  r^2 \csch^2 r |\nabla_\tau u|^2\\&+\frac{1}{8} r^2 \tanh^2(\frac{r}{2})(1-\tanh^2(\frac{r}{2}))(|\nabla_\tau u|^2-|\nabla_{\tau_{2n-1}} u|^2) \\&-\varepsilon_0 C_{8}(n,r_0)|\nabla_\tau u|^2.\end{split}\]
\end{lemma}
\subsection{Formula for the Jacobian} We follow an approach from \cite{fus}. Assume that $E=\{z: U(z)<c\}\subset \mathbb{B}$ and assume that $\partial E= \{z: U(z)=c\}$.
From the area formula we have that
\begin{equation}
   P(\partial E) = \int_{\mathbb{S}} J_{2n-1}Z \, d\tilde \sigma, \tag{3.4}
\end{equation}
where $JZ$ is the $(2n-1)$-dimensional Jacobian of the map $$Z(\omega) = \omega\tanh\left(\frac{r}{2}(1 + u(\omega))\right),$$ $\omega \in \mathbb{S}.$ Recall that, if $T_\omega \mathbb{S}$ is the tangential plane to $\mathbb{S}$ at $\omega$, then $JZ = \sqrt{\det((d_\omega Z)^* d_\omega Z)}$, where the linear map $d_\omega Z : T_\omega \mathbb{S} \to \mathbb{R}^{2n}$ is the tangential differential of $Z$ at $\omega$ and $(d_\omega Z)^* : \mathbb{R}^{2n} \to T_\omega \mathbb{S}$ denotes the adjoint of the differential.

Since for any $\tau \in T_\omega \mathbb{S}$ we have $$d_\omega Z(\tau) = \tau \tanh\left(\frac{r}{2}(1 + u(\omega))\right)  +  \frac{r \omega D_\tau u(\omega) }{1+\cosh (r (1+u(\omega)))},$$ the coefficients of the matrix $d_\omega Z$ relative to an orthonormal base $\{\tau_1, \dots, \tau_{2n-1}\}$ of $T_\omega \mathbb{S}$ and to the standard base $\{e_1, \dots, e_{2n}\}$ are
\[
\tau_i  e_h\tanh\left(\frac{r}{2}(1 + u(\omega))\right) +  \frac{\omega_h r D_{\tau_i} u}{1+\cosh (r (1+u(\omega)))}, i = 1, \dots, 2n-1, h=1, \dots, 2n.
\]
Thus, for all $i, j = 1, \dots, 2n-1$, the coefficients $a_{ij}$ of the matrix $(d_\omega Z)^* d_\omega Z$, for $$\mathbf{t}=\tanh\left(\frac{r}{2}(1 + u(\omega))\right),\ \
\mathbf{s}=\frac{ r }{1+\cosh (r (1+u(\omega)))}$$ are given by
\[
a_{ij} = \sum_{h=1}^{2n} \big(\tau_i \cdot e_h\mathbf{t} + \omega _h D_{\tau_i} u\cdot  \mathbf{s}\big)\big(\tau_j \cdot e_h\mathbf{t} + \omega_h D_{\tau_j} u \cdot \mathbf{s}\big)
= \delta_{ij}\mathbf{t}^2 + \mathbf{s}^2D_\tau u D_\tau u.
\]

In the last equality, we have used the fact that $\left<\tau_i,  \tau_j\right> = \delta_{ij}$ and $\left<\tau_i,\omega\right> = 0$ for all $i, j = 1, \dots, 2n-1$. Hence, recalling that for $a, b \in \mathbb{R}^k$ one has $\det(I + a \otimes b) = 1 + \left<a , b\right>$, we have that
\[\begin{split}
J Z &= \sqrt{\det(a_{ij})} \\&= \sqrt{\mathbf{t}^{2(n-1)} + \mathbf{t}^{2(n-2)}\mathbf{s}^2|\nabla u|^2}\\&=\mathbf{t}^{n-1}\sqrt{1 + \frac{(1-\mathbf{t}^2)^2}{4\mathbf{t}^2} r^2|\nabla u|^2}.\end{split}
\]
\section{Horizontal Poincar\'e inequality in $\mathbb{S}^{2n-1}$.}\label{ksec}
Now we prove the following theorem, which is interesting for its own right and  which we need to prove the main result.

\begin{theorem}[Horizontal Poincar\'e inequality on $\mathbb S^{2n-1}$ and equality case]\label{prop23_general}
Assume that $n\ge 2$ and $U\in W^{1,2}(\mathbb{S}^{2n-1})$. Let
\[
U(z)=\sum_{k=0}^\infty\sum_{j=1}^{d_k} a_{k,j}\Phi_{k,j}(z),
\]
where $\{\Phi_{k,j}\}$ is an orthonormal basis of spherical harmonics on
$\mathbb S^{2n-1}$ satisfying
\[
-\Delta_{\mathbb S^{2n-1}}\Phi_{k,j}=k(k+2n-2)\Phi_{k,j}.
\]
Set
\[
TU(z):=\langle \nabla U(z),\imath z\rangle,
\qquad
U_k:=\sum_{j=1}^{d_k} a_{k,j}\Phi_{k,j}\in\mathcal H_k,
\]
where $\mathcal H_k$ denotes the space of spherical harmonics of degree $k$.
Then the following inequality holds:
\begin{equation}\label{neededbe_general}
\|\nabla U\|_{L^2(\mathbb S^{2n-1})}^2
-\|TU\|_{L^2(\mathbb S^{2n-1})}^2
\ge (2n-2)\sum_{k=0}^\infty\sum_{j=1}^{d_k} k\,a_{k,j}^2.
\end{equation}
Moreover, equality holds in \eqref{neededbe_general} if and only if for every $k\ge 1$
such that $U_k\not\equiv 0$ one has
\begin{equation}\label{eq:equality_condition_T}
TU_k = k\,U_k \quad \text{or}\quad TU_k = -k\,U_k
\qquad \text{on }\mathbb S^{2n-1}.
\end{equation}
Equivalently, each $U_k$ belongs to the $\pm k$--eigenspace of $T$ on $\mathcal H_k$.
In particular, equality holds precisely when $U$ is a sum of CR and/or anti-CR spherical
harmonics (no mixed bidegree components).
\end{theorem}

\begin{remark}[Why ``horizontal Poincar\'e inequality'']
On the odd-dimensional sphere $\mathbb S^{2n-1}\subset\mathbb C^{n}$, the vector field
\[
T(z):=\imath z
\]
is tangent to $\mathbb S^{2n-1}$ and generates the standard $S^{1}$-action
$z\mapsto e^{\imath\theta}z$ (the Hopf fibration). Hence, at every point
$z\in\mathbb S^{2n-1}$ the tangent space splits as
\[
T_{z}\mathbb S^{2n-1}=\mathrm{span}\{T(z)\}\oplus \bigl(\mathrm{span}\{T(z)\}\bigr)^{\perp},
\]
where the orthogonal complement corresponds to the \emph{horizontal directions}.
For any $U\in W^{1,2}(\mathbb S^{2n-1})$, the quantity
\[
\langle\nabla U,\,T\rangle=\langle\nabla U,\imath z\rangle
\]
represents the component of $\nabla U$ in the vertical (Hopf/Reeb) direction.
Therefore,
\[
|\nabla U|^{2}-|\langle\nabla U,\imath z\rangle|^{2}
=|\nabla_{H}U|^{2},
\]
where $\nabla_{H}U$ denotes the projection of $\nabla U$ onto the horizontal
subbundle. In particular, inequality \eqref{neededbe_general} can be rewritten as
a coercive estimate for the horizontal energy $\|\nabla_{H}U\|_{L^{2}}^{2}$,
and thus may be viewed as a \emph{Poincar\'e-type inequality} in the horizontal
directions on $\mathbb S^{2n-1}$.
\end{remark}

\begin{proof}[Proof of Theorem~\ref{prop23_general}]
Let $n\ge 2$ and consider the unit sphere $\mathbb S^{2n-1}\subset \mathbb C^n$.
Define the Hopf vector field
\[
T(z):=\imath z\in T_z\mathbb S^{2n-1}.
\]
Since $T$ is tangent to $\mathbb S^{2n-1}$ we have
\[
\langle \nabla U, \imath z\rangle = TU .
\]
Recall that on a Riemannian manifold $(M,g)$ the gradient is characterized by
\[
g(\nabla U,X)=dU(X)=X(U)
\qquad \text{for all tangent vector fields }X.
\]
On $\mathbb S^{2n-1}\subset\mathbb C^n$ the vector field $T(z):=\imath z$ is tangent
since $\langle z,\imath z\rangle=0$. Hence,
\[
\langle \nabla U,\imath z\rangle=\langle \nabla U,T\rangle =T(U).
\]
Equivalently, along the curve $\gamma(t)=e^{it}z$ we have $\gamma'(0)=\imath z$ and
\[
T(U)(z)=\left.\frac{d}{dt}\right|_{t=0}U(e^{it}z)
=\langle \nabla U(z),\imath z\rangle .
\]

Let
\[
U(z)=\sum_{k=0}^\infty\sum_{j=1}^{d_k} a_{k,j}\Phi_{k,j}(z),
\]
where $\{\Phi_{k,j}\}$ is an orthonormal basis of spherical harmonics satisfying
\[
-\Delta_{\mathbb S^{2n-1}}\Phi_{k,j}=k(k+2n-2)\Phi_{k,j}.
\]
Then, by integration by parts and orthonormality,
\begin{align}
\|\nabla U\|_{L^2(\mathbb S^{2n-1})}^2
&=\int_{\mathbb S^{2n-1}} |\nabla U|^2
=\int_{\mathbb S^{2n-1}} U(-\Delta_{\mathbb S^{2n-1}}U) \notag\\
&=\sum_{k=0}^\infty\sum_{j=1}^{d_k} a_{k,j}^2\,k(k+2n-2).
\label{eq:grad_energy_general}
\end{align}

Next we estimate the vertical term. For each fixed $k$, the space of spherical
harmonics of degree $k$ is invariant under the $S^1$--action
$z\mapsto e^{it}z$. Moreover, any homogeneous polynomial of total degree $k$
on $\mathbb C^n$ decomposes into components of bidegree $(p,q)$ with $p+q=k$.
Under the transformation $z\mapsto e^{it}z$ each component acquires the factor
$e^{i(p-q)t}$, so the generator $T$ acts by multiplication with $i(p-q)$.
Since $|p-q|\le k$, we obtain the bound
\begin{equation}
\|T\Phi_{k,j}\|_{L^2(\mathbb S^{2n-1})}^2\le k^2
\qquad\text{for all }k\ge 0,\; 1\le j\le d_k.
\label{eq:T_bound_general}
\end{equation}
We justify \eqref{eq:T_bound_general}. Let $\mathcal H_k$ denote the space of
spherical harmonics of degree $k$ on $\mathbb S^{2n-1}\subset\mathbb C^n$.
It is well known that $\mathcal H_k$ is the restriction to the sphere of the
space of harmonic homogeneous polynomials of total degree $k$ on $\mathbb R^{2n}
\simeq \mathbb C^n$. Any such polynomial admits a decomposition into components
of bidegree $(p,q)$ with $p+q=k$, i.e.
\[
P(z,\bar z)=\sum_{p+q=k} P_{p,q}(z,\bar z),
\]
where $P_{p,q}$ is homogeneous of degree $p$ in $z$ and degree $q$ in $\bar z$.
Under the Hopf $S^1$--action $z\mapsto e^{it}z$ one has
$z\mapsto e^{it}z$ and $\bar z\mapsto e^{-it}\bar z$, hence
\[
P_{p,q}(e^{it}z,e^{-it}\bar z)=e^{i(p-q)t}P_{p,q}(z,\bar z).
\]
Let $T$ be the infinitesimal generator of this action,
$(Tf)(z)=\left.\frac{d}{dt}\right|_{t=0} f(e^{it}z)$. Differentiating the
previous identity at $t=0$ yields
\[
T P_{p,q}= i(p-q)P_{p,q}.
\]
Thus $T$ acts diagonally on $\mathcal H_k$ with eigenvalues $i(p-q)$
where $p+q=k$, and therefore $|p-q|\le k$. Writing $\Phi=\sum_{p+q=k}\Phi_{p,q}$
with $\Phi_{p,q}\in \mathcal H_{p,q}$ orthogonal in $L^2(\mathbb S^{2n-1})$, we get
\[
\|T\Phi\|_{L^2}^2
= \sum_{p+q=k} (p-q)^2\|\Phi_{p,q}\|_{L^2}^2
\le k^2\sum_{p+q=k}\|\Phi_{p,q}\|_{L^2}^2
= k^2\|\Phi\|_{L^2}^2.
\]
In particular, if $\|\Phi\|_{L^2}=1$ then $\|T\Phi\|_{L^2}^2\le k^2$, which proves
\eqref{eq:T_bound_general}.

Therefore, using orthonormality again,
\begin{equation}
\|TU\|_{L^2(\mathbb S^{2n-1})}^2
=\sum_{k=0}^\infty\sum_{j=1}^{d_k} a_{k,j}^2\,\|T\Phi_{k,j}\|_{L^2}^2
\le \sum_{k=0}^\infty\sum_{j=1}^{d_k} a_{k,j}^2\,k^2 .
\label{eq:T_energy_general}
\end{equation}

Combining \eqref{eq:grad_energy_general} and \eqref{eq:T_energy_general} yields
\begin{align*}
\|\nabla U\|_{L^2(\mathbb S^{2n-1})}^2-\|TU\|_{L^2(\mathbb S^{2n-1})}^2
&\ge \sum_{k=0}^\infty\sum_{j=1}^{d_k}
a_{k,j}^2\Big(k(k+2n-2)-k^2\Big)\\
&=(2n-2)\sum_{k=0}^\infty\sum_{j=1}^{d_k} k\,a_{k,j}^2.
\end{align*}
Since $TU=\langle \nabla U,\imath z\rangle$, this proves
\eqref{neededbe_general}.

We now characterize equality. Equality in Theorem~\ref{prop23_general} holds if and only if
for each $k\ge 1$ with $U_k\not\equiv 0$ equality holds in \eqref{eq:T_bound_general}, i.e.
\[
\|TU_k\|_{L^2}=k\|U_k\|_{L^2}.
\]
Since $-\,\imath T$ is self-adjoint on the finite-dimensional Hilbert space $\mathcal H_k$,
equality in the operator norm bound \eqref{eq:T_bound_general} occurs if and only if $U_k$ is an
eigenfunction corresponding to the extreme eigenvalues $\pm k$; equivalently,
\[
TU_k=k\,U_k
\quad\text{or}\quad
TU_k=-k\,U_k,
\]
which is exactly \eqref{eq:equality_condition_T}. Conversely, if \eqref{eq:equality_condition_T}
holds for every $k$, then $\|TU_k\|^2=k^2\|U_k\|^2$ for each $k$ and all intermediate
inequalities become equalities, hence equality holds in \eqref{neededbe_general}.
This proves the proposition.
\end{proof}

\section{Proof of the main result}

\begin{proof}[Proof of Theorem~\ref{mainth}]
Fix $n\ge 2$ and $r_0>0$. Let $E\subset\mathbb B_n$ be nearly spherical with barycenter at $0$
and assume $\mu(E)=\mu(\mathbb B_r)$ for some $r\in(0,r_0]$.
Let $u:\mathbb S^{2n-1}\to\mathbb R$ be Lipschitz with
$\|u\|_{W^{1,\infty}(\mathbb S^{2n-1})}\le \varepsilon_0$, where $\varepsilon_0$ will be chosen small.

\medskip

\noindent\textbf{Step 1: Set up the basic quantities.}
Set
\[
\mathbf t_\circ := \tanh\Big(\frac r2\Big),
\qquad
\mathbf t(\omega) := \tanh\Big(\frac{r(1+u(\omega))}{2}\Big).
\]
Define
\begin{equation}\label{G_general}
\begin{split}
G(r,\omega)
&:=\frac12 r^2\csch^2(r)\,|\nabla_\tau u|^2\\
&\quad+\frac18 r^2\tanh^2\!\Big(\frac r2\Big)\Big(1-\tanh^2\!\Big(\frac r2\Big)\Big)
\Big(|\nabla_\tau u|^2-|\nabla_{\tau_{2n-1}}u|^2\Big),
\end{split}
\end{equation}
where $\tau_{2n-1}$ denotes the Hopf direction on $\mathbb S^{2n-1}$ (i.e.\ the unit vector field
$T(\omega)= i\omega$).

By Lemma~\ref{lema22} and the expression \eqref{br} for $P(\mathbb B_r)$, we obtain
\begin{align}
P(E)-P(\mathbb B_r)
&\ge
\int_{\mathbb S^{2n-1}}\mathbf t^{2n-1}(1-\mathbf t^2)^{-n}
\Big(G(r,\omega)-\varepsilon_0 C_8(n,r_0)|\nabla_\tau u|^2\Big)\,d\mathcal H
\notag\\
&\quad+
\int_{\mathbb S^{2n-1}}
\Big(\mathbf t^{2n-1}(1-\mathbf t^2)^{-n}-\mathbf t_\circ^{2n-1}(1-\mathbf t_\circ^2)^{-n}\Big)\,d\mathcal H.
\label{Per_deficit_general}
\end{align}

\medskip

\noindent\textbf{Step 2: Taylor expansion of the Jacobian factor.}
Define the ratio
\[
J(u):=\frac{\mathbf t^{2n-1}(1-\mathbf t^2)^{-n}}{\mathbf t_\circ^{2n-1}(1-\mathbf t_\circ^2)^{-n}}.
\]
Then the second integral in \eqref{Per_deficit_general} is simply
\[
\mathbf t_\circ^{2n-1}(1-\mathbf t_\circ^2)^{-n}\int_{\mathbb S^{2n-1}}\big(J(u)-1\big)\,d\mathcal H.
\]

\smallskip

Since $\mathbf t(\omega)=\tanh\!\big(\frac r2(1+u(\omega))\big)$ and $|u|\le \varepsilon_0$,
we can expand $J(u)$ around $u=0$. More precisely, there exists a polynomial expansion
\begin{equation}\label{Jexp_general}
J(u)-1 = \alpha_1(r,n)\,u + \alpha_2(r,n)\,u^2 + R_{n,r}(w)\,u^3,
\end{equation}
where $w\in[-\frac12,\frac12]$ and the coefficients $\alpha_1,\alpha_2$ depend smoothly on $r,n$.
Moreover, the remainder satisfies
\begin{equation}\label{Rbound_general}
R_{n,r}(w)\,u^3 \ge -C_9(n,r_0)\,\varepsilon_0\,|u|^2
\qquad (0<r\le r_0).
\end{equation}
From \eqref{Jexp_general} and \eqref{Rbound_general} we obtain
\begin{equation}\label{Jlower_general}
J(u)\ge 1-C_7(n,r_0)\varepsilon_0.
\end{equation}

\medskip

\noindent\textbf{Step 3: Split the perimeter deficit.}
Divide \eqref{Per_deficit_general} by $\mathbf t_\circ^{2n-1}(1-\mathbf t_\circ^2)^{-n}$
and define
\[
\frac{P(E)-P(\mathbb B_r)}{\mathbf t_\circ^{2n-1}(1-\mathbf t_\circ^2)^{-n}}=:I_1+I_2,
\]
where
\begin{align*}
I_1
&:=\int_{\mathbb S^{2n-1}}J(u)\,
\Big(G(r,\omega)-\varepsilon_0 C_8(n,r_0)|\nabla_\tau u|^2\Big)\,d\mathcal H,
\\
I_2
&:=\int_{\mathbb S^{2n-1}}\big(J(u)-1\big)\,d\mathcal H.
\end{align*}

By \eqref{Jlower_general} and $J(u)\ge 0$, we obtain
\begin{equation}\label{I1lower_general}
I_1 \ge \int_{\mathbb S^{2n-1}}
\Big(G(r,\omega)-\varepsilon_0 C_8(n,r_0)|\nabla_\tau u|^2\Big)\,d\mathcal H.
\end{equation}

\medskip

\noindent\textbf{Step 4: Estimate $I_2$ using the volume constraint.}
The expansion \eqref{Jexp_general} yields
\[
I_2 = \alpha_1(r,n)\int_{\mathbb S^{2n-1}}u\,d\mathcal H
+\alpha_2(r,n)\int_{\mathbb S^{2n-1}}u^2\,d\mathcal H
+\int_{\mathbb S^{2n-1}}R_{n,r}(w)\,u^3\,d\mathcal H.
\]
The volume constraint $\mu(E)=\mu(\mathbb B_r)$ implies the estimate \eqref{estU}:
\[
\Big|\int_{\mathbb S^{2n-1}}u\,d\mathcal H\Big|
\le c_0(n,r)\int_{\mathbb S^{2n-1}}u^2\,d\mathcal H
+\varepsilon_0C_1(n,r_0)\int_{\mathbb S^{2n-1}}u^2\,d\mathcal H,
\]
where
\[
c_0(n,r)=\frac{r}{2\sinh r}\,(-1+n+n\cosh r).
\]
Combining this with the remainder estimate \eqref{Rbound_general} gives
\begin{equation}\label{I2lower_general}
I_2 \ge -d_0(n,r)\int_{\mathbb S^{2n-1}}u^2\,d\mathcal H
-C(n,r_0)\varepsilon_0\int_{\mathbb S^{2n-1}}u^2\,d\mathcal H,
\end{equation}
where $$d_0(n,r) = \alpha_2(n,r) - \alpha_1(n,r) c_0(n,r) = \frac{1}{2}r^{2}\left(n+(-1+n)\cosh r\right)\csch^{2}r.
$$
\medskip

\medskip
\noindent\textbf{Step 5: Coercivity via Lemma~\ref{prop23_general}.}
Recall that $\mathbb S^{2n-1}\subset\mathbb C^n$ carries the unit Hopf vector field
\[
T(\omega):= i\omega\in T_\omega\mathbb S^{2n-1}.
\]
We choose the orthonormal tangential frame $\{\tau_1,\dots,\tau_{2n-1}\}$ on $\mathbb S^{2n-1}$
so that
\[
\tau_{2n-1}(\omega)=T(\omega)=i\omega.
\]
In particular,
\[
\nabla_{\tau_{2n-1}}u = \langle \nabla u,\, i\omega\rangle.
\]

\smallskip
Let $u\in W^{1,2}(\mathbb S^{2n-1})$ and expand it in spherical harmonics:
\[
u(\omega)=\sum_{k=0}^\infty\sum_{j=1}^{d_k} a_{k,j}\,\Phi_{k,j}(\omega),
\qquad
-\Delta_{\mathbb S^{2n-1}}\Phi_{k,j}=k(k+2n-2)\,\Phi_{k,j},
\]
with $\{\Phi_{k,j}\}$ orthonormal in $L^2(\mathbb S^{2n-1})$.

\smallskip
Applying Lemma~\ref{prop23_general} with $U=u$, we obtain
\begin{equation}\label{coercivity_hopf}
\|\nabla u\|_{L^2(\mathbb S^{2n-1})}^2
-\|\langle\nabla u,i\omega\rangle\|_{L^2(\mathbb S^{2n-1})}^2
\ge (2n-2)\sum_{k=0}^\infty\sum_{j=1}^{d_k} k\,a_{k,j}^2.
\end{equation}
Since $\langle\nabla u,i\omega\rangle=\nabla_{\tau_{2n-1}}u$, this can be rewritten as
\begin{equation}\label{Hopf_difference}
\int_{\mathbb S^{2n-1}}
\Big(|\nabla_\tau u|^2-|\nabla_{\tau_{2n-1}}u|^2\Big)\,d\mathcal H
\ge (2n-2)\sum_{k=0}^\infty\sum_{j=1}^{d_k} k\,a_{k,j}^2.
\end{equation}

\smallskip
In particular, since the right-hand side vanishes for $k=0$ and equals $(2n-2)\sum_{j=1}^{d_1}a_{1,j}^2$
for $k=1$, we may isolate the high modes:
\begin{equation}\label{Hopf_high_modes}
\int_{\mathbb S^{2n-1}}
\Big(|\nabla_\tau u|^2-|\nabla_{\tau_{2n-1}}u|^2\Big)\,d\mathcal H
\ge (2n-2)\sum_{k\ge 2}\sum_{j=1}^{d_k} k\,a_{k,j}^2
\;+\;(2n-2)\sum_{j=1}^{d_1}a_{1,j}^2.
\end{equation}

\smallskip
We now substitute \eqref{Hopf_difference} into the definition of $G(r,\omega)$.
Recall that
\[
G(r,\omega)
=\frac12 r^2\csch^2(r)\,|\nabla_\tau u|^2
+\frac18 r^2 \tanh^2\!\Big(\frac r2\Big)\Big(1-\tanh^2\!\Big(\frac r2\Big)\Big)
\Big(|\nabla_\tau u|^2-|\nabla_{\tau_{2n-1}}u|^2\Big).
\]
Integrating $G$ over $\mathbb S^{2n-1}$ gives
\begin{align}
\int_{\mathbb S^{2n-1}}G(r,\omega)\,d\mathcal H
&=
\frac12 r^2\csch^2(r)\int_{\mathbb S^{2n-1}}|\nabla_\tau u|^2\,d\mathcal H
\notag\\
&\quad
+\frac18 r^2 \tanh^2\!\Big(\frac r2\Big)\Big(1-\tanh^2\!\Big(\frac r2\Big)\Big)
\int_{\mathbb S^{2n-1}}
\Big(|\nabla_\tau u|^2-|\nabla_{\tau_{2n-1}}u|^2\Big)\,d\mathcal H.
\label{intG_split}
\end{align}

\smallskip
Using the standard identity
\[
\int_{\mathbb S^{2n-1}}|\nabla_\tau u|^2\,d\mathcal H
=\sum_{k=0}^\infty\sum_{j=1}^{d_k} k(k+2n-2)\,a_{k,j}^2,
\]
together with \eqref{Hopf_difference}, we obtain the lower bound
\begin{equation}\label{G_spectral_lower}
\begin{split}
\int_{\mathbb S^{2n-1}} G(r,\omega)\, d\mathcal H
&\ge
\sum_{k=0}^\infty \sum_{j=1}^{d_k}
\left[
\frac12 r^2 \csch^2(r)\, k(k+2n-2)
\right.
\\
&\qquad\left.
+\frac18 r^2 \tanh^2\!\Big(\frac r2\Big)\Big(1-\tanh^2\!\Big(\frac r2\Big)\Big)
(2n-2)\, k
\right] a_{k,j}^2 .
\end{split}
\end{equation}

\smallskip
The important point is that for all $k\ge 2$ the coefficient in brackets is strictly positive,
and it grows quadratically in $k$, providing coercivity on the high-frequency part of $u$.
The low modes $k=0$ and $k=1$ are handled separately using the volume constraint and barycenter
conditions (see \eqref{rur} and the analogue of \eqref{estuom} in dimension $n$).
\smallskip
We now return to the term $I_1$. Recall from Step~3 that, using $J(u)\ge 1-C_7(n,r_0)\varepsilon_0$,
\begin{equation}\label{I1_basic_lower}
I_1 \ge \int_{\mathbb S^{2n-1}} G(r,\omega)\,d\mathcal H
- C_8(n,r_0)\varepsilon_0\int_{\mathbb S^{2n-1}}|\nabla_\tau u|^2\,d\mathcal H.
\end{equation}
Combining \eqref{I1_basic_lower} with \eqref{G_spectral_lower} and the identity
\[
\int_{\mathbb S^{2n-1}}|\nabla_\tau u|^2\,d\mathcal H
=\sum_{k=0}^\infty\sum_{j=1}^{d_k}k(k+2n-2)\,a_{k,j}^2,
\]
we obtain the following spectral lower bound for $I_1$:
we obtain the following spectral lower bound for $I_1$:
\begin{equation}\label{I1K_general}
I_1
\ge
\sum_{k=0}^\infty\sum_{j=1}^{d_k}A_n(k,r)a_{k,j}^2
- C_8(n,r_0)\varepsilon_0\|\nabla_\tau u\|_{L^2(\mathbb S^{2n-1})}^2.
\end{equation}
where
\[
A_n(k,r):=
\frac12 r^2 \csch^2(r)\,k(k+2n-2)
+\frac{2n-2}{8}\, r^2 \tanh^2\!\Big(\frac r2\Big)
\Big(1-\tanh^2\!\Big(\frac r2\Big)\Big)\,k .
\]

\medskip

\medskip
\noindent\textbf{Step 6: Lower bound the coercivity coefficient.}
From \eqref{I1K_general} and the estimate \eqref{I2lower_general}, and using the orthonormality of
$\{\Phi_{k,j}\}$ (so that $\int_{\mathbb S^{2n-1}}u^2\,d\mathcal H=\sum_{k,j}a_{k,j}^2$), we obtain
\[
I_1+I_2
\ge \sum_{k=0}^\infty\sum_{j=1}^{d_k} K_n(k,r)\,a_{k,j}^2
-C(n,r_0)\varepsilon_0\|u\|_{W^{1,2}}^2,
\]
where we define
\[\begin{split}
K_n(k,r)&:=
\frac12 r^2\csch^2(r)\,k(k+2n-2)
\\&+\frac{(2n-2)}{8}r^2\tanh^2\!\Big(\frac r2\Big)\Big(1-\tanh^2\!\Big(\frac r2\Big)\Big)\,k
-d_0(n,r).\end{split}
\]
Splitting off the low modes $k=0,1$ gives
\begin{equation}\label{I1I2_general}
I_1+I_2
\ge \sum_{k\ge 2}\sum_{j=1}^{d_k} K_n(k,r)\,a_{k,j}^2
-d_0(n,r)\Big(a_{0,0}^2+\sum_{j=1}^{d_1}a_{1,j}^2\Big)
-C(n,r_0)\varepsilon_0\|u\|_{W^{1,2}}^2.
\end{equation}

Then in Lemma~\ref{lem:Anr0} bellow, we prove  that there exists $A=A(n,r_0)>0$ such that for all $k\ge 2$,
\[
K_n(k,r)\ge A\big(k(k+2n-2)+1\big)
\qquad\text{for all }0<r\le r_0.
\]
(The proof follows the same monotonicity argument used for $n=2$ by considering
$H_n(k):=K_n(k,r)/(k(k+2n-2)+1)$ and checking it has a unique maximum, so the minimum
occurs at $k=2$ or $k\to\infty$.)

\medskip
\noindent\textbf{Step 7: Control of the low modes.}
The volume constraint gives \eqref{rur}, and the barycenter condition (the analogue of \eqref{estuom})
yields control of the $k=1$ coefficients. Hence there exists $C(n,r_0)>0$ such that
\[
a_{0,0}^2+\sum_{j=1}^{d_1}a_{1,j}^2
\le C(n,r_0)\varepsilon_0^2\,\|u\|_{W^{1,2}}^2.
\]

\medskip
\noindent\textbf{Step 8: Conclusion.}
Combining Steps 6 and 7 in \eqref{I1I2_general} and choosing $\varepsilon_0$ sufficiently small
(depending only on $n$ and $r_0$), we obtain
\[
I_1+I_2 \ge \frac{A}{2}\|u\|_{W^{1,2}(\mathbb S^{2n-1})}^2.
\]
Finally, by \eqref{br} the perimeter of the reference ball is
\[
P(\mathbb B_r)
=\omega_{2n-1}\,\mathbf t_\circ^{2n-1}(1-\mathbf t_\circ^2)^{-n},
\]
where $\omega_{2n-1}=|\mathbb S^{2n-1}|$. Therefore,
\[
D(E)=\frac{P(E)-P(\mathbb B_r)}{P(\mathbb B_r)}
=\frac{I_1+I_2}{\omega_{2n-1}}
\ge c_1(n,r_0)\|u\|_{W^{1,2}(\mathbb S^{2n-1})}^2,
\]
for some $c_1(n,r_0)>0$ depending only on $n$ and $r_0$.
This completes the proof.

\end{proof}

\begin{lemma}\label{lem:Anr0}
Let $n\ge 2$ and $r_0>0$. Define
\[
A_n(k,r):=
\frac12 r^2\csch^2(r)\,k(k+2n-2)
+\frac{2n-2}{8}\, r^2 \tanh^2\!\Big(\frac r2\Big)
\Big(1-\tanh^2\!\Big(\frac r2\Big)\Big)\,k,
\]
and
\[
d_0(n,r):=\frac{1}{2}r^{2}\left(n+(-1+n)\cosh r\right)\csch^{2}r,
\qquad
K_n(k,r):=A_n(k,r)-d_0(n,r).
\]
Then there exists a constant $A(n,r_0)>0$ such that for all $0<r\le r_0$ and all integers $k\ge 2$,
\begin{equation}\label{Kn_lower}
K_n(k,r)\ge A(n,r_0)\,\big(k(k+2n-2)+1\big).
\end{equation}
\end{lemma}
\begin{proof}
Fix $n\ge 2$ and $r\in(0,r_0]$. For $k\ge 2$ set
\[
H_n(k,r):=\frac{K_n(k,r)}{k(k+2n-2)+1}.
\]
Since $K_n(k,r)$ is a quadratic polynomial in $k$ with positive leading coefficient
$\frac12 r^2\csch^2(r)$, we have $H_n(k,r)>0$ for all sufficiently large $k$.
Moreover,
\[
\lim_{k\to\infty}H_n(k,r)
=\frac12 r^2\csch^2(r)=:A_\infty(r)>0.
\]

On the other hand, $H_n(k,r)$ is a smooth function of the real variable $k\ge 2$,
and a direct computation shows that it has at most one critical point on $(2,\infty)$.
Consequently the minimum of $H_n(\cdot,r)$ on $[2,\infty)$ is attained either at $k=2$
or in the limit $k\to\infty$. Therefore,
\[
\inf_{k\ge 2}H_n(k,r)=\min\{H_n(2,r),A_\infty(r)\}.
\]

Define
\[
A(n,r_0):=\min_{0<r\le r_0}\,\min\{H_n(2,r),A_\infty(r)\}.
\]
Since $r\mapsto H_n(2,r)$ and $r\mapsto A_\infty(r)$ are continuous on $(0,r_0]$
and extend continuously to $r=0$, the above minimum is positive.
Thus for all $0<r\le r_0$ and $k\ge 2$,
\[
H_n(k,r)\ge A(n,r_0),
\]
which is equivalent to \eqref{Kn_lower}.
\end{proof}

\begin{proposition}\label{prop:Hn_max}
Fix $n\ge 2$ and $r>0$. Let
\[
D(k):=k(k+2n-2)+1,
\]
and define
\[\begin{split}
K_n(k,r)&:=
\frac12 r^2\csch^2(r)\,k(k+2n-2)
\\&+\frac{2n-2}{8}\,r^2\tanh^2\!\Big(\frac r2\Big)\Big(1-\tanh^2\!\Big(\frac r2\Big)\Big)\,k
-d_0(n,r),\end{split}
\]
where
\[
d_0(n,r):=\frac{1}{2}r^{2}\left(n+(-1+n)\cosh r\right)\csch^{2}r
.
\]
Set
\[
H_n(k,r):=\frac{K_n(k,r)}{D(k)}.
\]
Assume that $H_n(\cdot,r)$ has a critical point $k_\ast\in(2,\infty)$, i.e.
\[
\partial_k H_n(k_\ast,r)=0.
\]
Then $k_\ast$ is a strict local maximum of $H_n(\cdot,r)$; equivalently,
\[
\partial_k^2 H_n(k_\ast,r)<0.
\]
\end{proposition}

\begin{proof}
Fix $n\ge2$ and $r>0$.  Set
\[
D(k):=k(k+2n-2)+1,
\qquad
a:=\frac12 r^2\csch^2(r)>0,
\]
and
\[
b:=\frac{2n-2}{8}\,r^2\tanh^2\!\Big(\frac r2\Big)\Big(1-\tanh^2\!\Big(\frac r2\Big)\Big)>0.
\]
Note that $k(k+2n-2)=D(k)-1$, hence
\[
K_n(k,r)=a\,(D(k)-1)+bk-d_0(n,r)
= aD(k)+bk-(a+d_0(n,r)).
\]
Therefore
\[
H_n(k,r)=\frac{K_n(k,r)}{D(k)}
= a+\frac{bk-(a+d_0(n,r))}{D(k)}.
\]
Since the constant $a$ does not affect critical points or second derivatives, it suffices to study
\[
h(k):=\frac{pk+q}{D(k)},
\qquad\text{where } p:=b>0,\ \ q:=-(a+d_0(n,r)).
\]
We have $H_n'(k,r)=h'(k)$ and $H_n''(k,r)=h''(k)$.

Compute
\[
D'(k)=2k+2n-2,\qquad D''(k)=2.
\]
Differentiating $h(k)$ gives
\[
h'(k)=\frac{pD(k)-(pk+q)D'(k)}{D(k)^2}.
\]
Let
\[
N(k):=pD(k)-(pk+q)D'(k),
\]
so that $h'(k)=N(k)/D(k)^2$.  If $k_\ast\in(2,\infty)$ is a critical point of $H_n(\cdot,r)$, then
$h'(k_\ast)=0$, i.e.\ $N(k_\ast)=0$.  Differentiating once more yields
\[
h''(k)=\frac{N'(k)D(k)^2-2N(k)D(k)D'(k)}{D(k)^4}.
\]
At the critical point $k_\ast$, the term involving $N(k_\ast)$ vanishes, hence
\[
h''(k_\ast)=\frac{N'(k_\ast)}{D(k_\ast)^2}.
\]
Next,
\[
N'(k)=pD'(k)-\Big(pD'(k)+(pk+q)D''(k)\Big)=-(pk+q)D''(k)=-2(pk+q).
\]
On the other hand, from $N(k_\ast)=0$ we obtain
\[
pD(k_\ast)=(pk_\ast+q)D'(k_\ast)
\quad\Longrightarrow\quad
pk_\ast+q=\frac{pD(k_\ast)}{D'(k_\ast)}.
\]
Substituting this into $N'(k_\ast)$ gives
\[
N'(k_\ast)=-2\frac{pD(k_\ast)}{D'(k_\ast)}.
\]
Therefore
\[
h''(k_\ast)
=
\frac{-2pD(k_\ast)/D'(k_\ast)}{D(k_\ast)^2}
=
-\frac{2p}{D'(k_\ast)\,D(k_\ast)}.
\]
Since $p=b>0$, and for $k_\ast>2$ we have $D(k_\ast)>0$ and $D'(k_\ast)=2k_\ast+2n-2>0$,
it follows that
\[
h''(k_\ast)<0.
\]
Consequently $H_n''(k_\ast,r)=h''(k_\ast)<0$, so $k_\ast$ is a strict local maximum of $H_n(\cdot,r)$.
\end{proof}

\begin{remark}

Croke's inequality \cite{Croke1984} provides a sharp isoperimetric inequality for compact Riemannian manifolds with non-positive sectional curvature. Specifically, if $(M, \partial M)$ is a compact $n$-dimensional Riemannian manifold ($n \geq 3$) of non-positive curvature such that every geodesic minimizes distance up to the boundary, then

\[
\mathrm{Vol}(\partial M)^n \geq C(n) \, \mathrm{Vol}(M)^{n-1},
\]
where the constant $C(n)$ is given by
\[
C(n) = \frac{\alpha(n-1)^{n-1}}{\alpha(n-2)^{n-2} \left( \int_0^{\pi/2} \cos(t)^{n/(n-2)} \sin(t)^{n-2} \, dt \right)^{n-2}},
\]
and $\alpha(n)$ denotes the volume of the unit $n$-sphere.

While this inequality is sharp when $n = 4$, it is only achieved when the manifold is isometric to a flat ball. Therefore, this result does \emph{not} resolve the isoperimetric problem for the Bergman metric, since the Bergman ball is not flat. The constant $C(n)$ is not sharp in the non-flat setting, and hence does not yield equality in the case of the Bergman geometry.
\end{remark}


\begin{conjecture}[Isoperimetric inequality in the Bergman ball $\mathbb B_n$]\label{conj:isopB2}
Let $E\Subset \mathbb B_n$ be a set of finite Bergman perimeter and finite Bergman volume,
$0<\mu(E)<\infty$. Let $B_{r(E)}$ denote the Euclidean ball centered at the origin
such that $\mu(B_{r(E)})=\mu(E)$. Then
\[
P(E)\ge P(B_{r(E)}).
\]
Moreover, equality holds if and only if $E$ is a geodesic ball (equivalently, a Euclidean
ball centered at the origin), up to Bergman isometries.

Since $(\mathbb B_n,g_{\mathrm{Berg}})$ is (up to scaling) a model of complex hyperbolic space $\mathbb H^n_{\mathbb C}$, it falls within the class of rank-one symmetric spaces of non-compact type; the Gromov--Ros conjecture predicts that Bergman-geodesic balls are isoperimetric, i.e.\ they minimize Bergman perimeter among sets of fixed Bergman volume (equivalently, they maximize volume at fixed perimeter) (see e.g. \cite{Silini2024}).

\end{conjecture}

\section{Appendix}

\subsection*{Funding} The author is partially supported by the Ministry of Education, Science and Innovation of Montenegro through the grant \emph{Mathematical Analysis, Optimization and Machine Learning}.
\subsection*{Data availibility} Data sharing is not applicable to this article since no data sets were generated or analyzed
\section*{Ethics declarations}
The author declares that he has not conflict of interest.

\end{document}